\documentclass{article} 
\usepackage{amsmath,amsfonts,amsbsy,amsthm,latexsym,amssymb,enumerate}

\newtheorem{theorem}{Theorem}[section] 
\newtheorem{corollary}[theorem]{Corollary} 
\newtheorem{lemma}[theorem]{Lemma} 
\newtheorem{proposition}[theorem]{Proposition}

\begin{document} 
\sloppy 
\title{~\\[-6ex] Limit theorems for weakly subcritical branching  
processes in random environment\thanks{This paper is part of a  project 
supported by the German Research Foundation (DFG) and the 
Russian Foundation of Basic Research (Grant DFG-RFBR 08-01-91954)}} 
\author{\textsc{V.I. Afanasyev}\thanks{Department of 
Discrete Mathematics, Steklov Mathematical Institute, 8 Gubkin Street, 
119\,991 Moscow, Russia, viafan@mail.ru} 
   \hspace{.8cm} 
\textsc{C. B\"oinghoff} 
 \thanks{Fachbereich Mathematik, Universit\"at Frankfurt, Fach 
    187, D-60054 Frankfurt am Main, Germany , boeinghoff@math.uni-frankfurt.de} \\ 
  \textsc{G. Kersting} 
   \thanks{Fachbereich Mathematik, Universit\"at Frankfurt, Fach 
    187, D-60054 Frankfurt am Main, Germany, 
 kersting@math.uni-frankfurt.de} 
     \hspace{.8cm} \textsc{V.A. Vatutin}\thanks{Department of 
Discrete Mathematics, Steklov Mathematical Institute, 8 Gubkin Street, 
119\,991 Moscow, Russia,  vatutin@mi.ras.ru }} 
\maketitle 
\begin{abstract} 
For a branching process in random environment it is assumed that the offspring distribution of the individuals varies in a random fashion, independently from one generation to the other. Interestingly there is the possibility that the process may at the same time be subcritical and, conditioned on nonextinction, 'supercritical'. This so-called weakly subcritical case is considered in this paper. We study the asymptotic survival probability and the size of the population conditioned on non-extinction. Also a functional limit theorem is proven, which makes the conditional supercriticality manifest. A main tool is a new type of functional limit theorem for conditional random walks.
\end{abstract} 
 \begin{small} 
 \emph{MSC 2000 subject classifications.}  Primary 60J80, Secondary 
60G50, 60F17  \\ 
\emph{Key words and phrases.} Branching process, random environment, 
 random walk, change of measure, survival probability, functional limit theorem
\end{small} 
\thispagestyle{empty} \newpage
\section{Introduction and main results} 
For a branching process 
in random environment 
it is assumed that the offspring distribution of the individuals varies in a random fashion, independently from one generation to the other. 
Conditioned  on 
the environment       individuals reproduce independently of  each 
other. Let $Q_{n}$ be the random offspring distribution 
 of an individual at generation~$n-1$ and let 
$Z_n$ denote the number 
of individuals at generation $n$. Then $Z_{n}$ is the sum of 
$Z_{n-1}$ independent random variables, each of which has 
distribution~$Q_{n}$. To give a formal definition let $\Delta$ 
be the 
   space of      probability measures on 
 \mbox{$\mathbb{N}_{0}:= \{0,1,\ldots\}$}, which equipped with the metric of     total   variation is 
a Polish space. Let $Q$ be a random 
   variable taking values in $\Delta$. 
Then, an infinite sequence $\Pi=(Q_{1}, 
  Q_{2},\ldots)$ of i.i.d.\ copies of  $Q$ is said to form a 
\emph{random environment}. A sequence of $\mathbb{N}_0$-valued random 
   variables      $Z_{0},Z_{1},\ldots$ is called a 
\emph{branching process in the random environment} $\Pi$, if $Z_{0} 
 $ is      independent of $ \Pi$   and   given 
  $ \Pi$ the  process $Z=(Z_{0},Z_{1},\ldots)$ is a Markov chain 
with 
\begin{equation}  \label{transition} 
   \mathcal{L} \big(Z_{n} \; \big| \; Z_{n-1}=z, \, \Pi 
  = 
    (q_{1},q_{2},\ldots) \big) \ = \ q_{n}^{*z} 
\end{equation} 
 for every $n\in \mathbb{N}$, $z \in \mathbb{N}_0$ and 
    $q_{1},q_{2},\ldots \in \Delta$, where 
$q^{*z}$  is the $z$-fold 
convolution of the measure $q$. 
The   corresponding probability 
   measure on the 
   underlying probability space will be denoted by $\mathbb{P}$. 
In the following 
we assume  that the process starts with a single founding 
ancestor, $Z_0=1$ a.s.,  and (without loss of generality) that 
$\mathbb{P}\{Q(\{0\})=1\}  =0 .$  Note that in general~$Z$ is not the superposition of $Z_0$ independent copies of the  process started at $Z_0=1$. 

It turns out  that the asymptotic behavior 
of the generation size  process $Z$ is  determined  in the main 
by the associated 
 random walk $S=(S_n)_{n\ge 0}$. This random walk has 
 initial state 
$S_{0}=0$ and  increments   $X_n=S_n-S_{n-1}, \, n\ge 1$ 
defined as 
\[ 
   X_{n} \ := \ \log m(Q_n),  
\] 
where 
\[ 
m(q)\ := \
\sum_{y=0}^{\infty} y  q(\{y\}) 
\] 
is the mean of the offspring distribution $q \in \Delta$. 
In view of (\ref{transition}) and the assumption $Z_0=1$ a.s.\ 
     the conditional expectation of 
$Z_{n}$      given the environment~$ \Pi$ 
  can be expressed by means of $(S)_{n\in\mathbb{N}_0}$ as 
\begin{equation}  \label{expect} 
 \mathbb{E}[Z_{n} \,| \,  \Pi \,]   \ = \ 
\prod_{k=1}^n m(Q_k)\ = \ \exp(S_n) \quad 
\mathbb{P}\text{--}a.s. 
\end{equation} 
Averaging over the environment gives 
\begin{equation}  \label{expect2} 
 \mathbb{E}[Z_{n} ]   \ = \ 
 \big(\mathbb{E}[ m(Q) ]\big)^n.
\end{equation} 

 If the random walk $(S)_{n\in\mathbb{N}_0}$ drifts to $-\infty$,  then the branching process 
is said to be   \emph{subcritical}. In case $X=\log  m(Q)$ has finite 
mean, subcriticality corresponds to $\mathbb{E}[X]<0.$ 
For such processes 
 the conditional non-extinction probability at $n$ 
 decays at an exponential rate for almost every environment. 
This fact is an immediate consequence of the strong law of large numbers and the first moment  
estimate 
\begin{eqnarray} 
  \lefteqn{ \mathbb{P} \{Z_{n} > 0 \,| \,  \Pi \} \ =  \ 
\min_{0\le k \leq n} 
      \mathbb{P} \{Z_{k} > 0 \, | \,  \Pi \}}     \nonumber\\ 
   &\leq & 
 \min_{0\le k \leq n}   \mathbb{E}[ Z_{k} \,| \,  \Pi ]  \ = \ 
\exp\big(\min_{0\le k \leq n} S_{k}\big)  \quad 
\mathbb{P}\text{--}a.s.       \label{gleichung1} 
\end{eqnarray} 

As was  first observed by Afanasyev~\cite{af_80} 
  and later independently by Dekking~\cite{de} there are three possibilities for the asymptotic behavior of subcritical branching processes. They are called the \emph{weakly} subcritical, the \emph{intermediate}  subcritical and the 
\emph{strongly} subcritical case.  Here we study the weakly subcritical case. 

The present article is part of several publications
having started  with \cite{agkv,agkv2}, in which we try to develop characteristic properties of the different  cases. For a comparative discussion 
we refer the reader to \cite{bgk}.
One purpose of this paper is to make the methods, developed in \cite{agkv} for criticality, also available for weak subcriticality.

\paragraph{Assumption A1. } {\em The process $Z$ is weakly subcritical, that is, there is a number $0 < \beta < 1$ such that}
 \[  \mathbb{E} [ Xe^{\beta X}] \ = \ 0 \ . \]  \smallskip

This implies $-\infty\leq\mathbb{E}[X]<0$, thus $(S)_{n\in\mathbb{N}_0}$ has negative drift with respect to $\mathbb{P}$. Assumption A1 is somewhat weaker than $\mathbb{E}[X] < 0 < \mathbb{E}[Xe^X] < \infty$, which is a customary condition for weak subcriticality. The assumption suggests to change to the measure $\mathbf{P}$ with expectation $\mathbf{E}$. For any $n\in\mathbb{N}$ and any measurable, bounded function $\varphi:\Delta^n\times \mathbb{N}_0^{n+1}\rightarrow\mathbb{R}$, the measure $\mathbf{P}$ is given by
\[ \mathbf{E}[ \varphi(Q_1,\ldots,Q_n, Z_0,\ldots,Z_n)] \ := \ \gamma^{-n} \mathbb{E}\big[\varphi(Q_1,\ldots,Q_n, Z_0,\ldots,Z_n)e^{\beta (S_n-S_0)}\big] \ , \] 
with
\[ \gamma \ := \ \mathbb{E} [e^{\beta X}] \ . \]  
(We include $S_0$ in the above expression, because later we shall also consider cases where $S_0 \neq 0$.) Then $ \mathbb{E}[Xe^{\beta X}]=0$ translates into
\[ \mathbf{E}[X] \ = \ 0 \ . \]
Thus $(S)_{n\in\mathbb{N}_0}$ becomes a recurrent random walk under $\mathbf{P}$. 

As to the regularity of the distribution of $X$ we make the following assumptions.

\paragraph{Assumption A2.} {\em The distribution of $X$ has  finite variance with respect to $\mathbf{P}$ or  (more generally) belongs to the domain of attraction of some stable law with index $\alpha \in (1,2]$. It is non-lattice.} \\

\noindent
Since $\mathbf{E}[X]=0$ this means that there is an increasing sequence of positive numbers
\[ a_n  \ = \ n^{1/\alpha} \ell_n \]
with a slowly varying sequence $\ell_1,\ell_2,\ldots$ such that for $n\rightarrow\infty$
\[ \mathbf{P}\{S_n/a_n \in dx \, \} \ \to \ s(x) \, dx\] 
weakly, where $s(x)$ denotes the density of the limiting stable law. Below the local version of this statement will mainly be used. Note that due to the change of measure $X^-$ always has finite variance and an infinite variance may only arise from $X^+$. Then the stable law $s(x) \, dx$ is completely skewed towards the positive real axis. Nevertheless it is not a one sided law: $s(x)$ is strictly positive everywhere, as this is the case for stable laws with finite expectation. 

\paragraph{Remark.} In \cite{agkv} we studied branching processes in a critical random environment under the assumption that the random walk $(S)_{n\in\mathbb{N}_0}$ fulfils {\em Spitzer's condition}. In general this condition is less restrictive than A2. However, if $X^-$ has finite second moment, then A2 is equivalent to Spitzer's condition (cf. \cite{do}). \hfill $\Box$ \\ 

Our last assumption on the environment concerns the standardized truncated second moment of $Q$,
\[ \zeta(a) \ := \ \sum_{y=a}^\infty y^2 Q(\{y\}) \Big/ m(Q)^2 \ , \quad a \in \mathbb{N} \ . \]

\paragraph{Assumption A3.} {\em For some $\varepsilon > 0$ and some $a \in \mathbb{N}$}
\[ \mathbf{E} [ (\log^+ \zeta(a))^{\alpha + \varepsilon}] \ < \ \infty \ ,\]
where $\log^+ x:= \log(\max(x,1))$.
\paragraph{Remark} For examples where this assumption is fulfilled, see \cite{agkv}. In particular our results hold for binary branching processes in random environment (where individuals have either two children or none) and for cases where $Q$ is a.s. a Poisson distribution or a.s. a geometric distribution.\\

We now come to the main results of the paper. All our limit theorems are under the law $\mathbb{P}$ which is what is called the annealed approach. The first theorem describes the asymptotic behaviour of the non--extinction probability at generation $n$. In the following, for sequences $(d_n)$ and $(m_n)$, we write $d_n\sim m_n$ if $d_n/m_n\rightarrow 1$ as $n\rightarrow\infty$.

\begin{theorem} \label{th1} Assume A1 to A3. Then there exists a number $0 < \kappa < \infty$ such that
\[ \mathbb{P} \{Z_n > 0\} \ \sim \ \kappa \ \mathbb{P} \{ \min(S_1,\ldots,S_n) \ge 0 \} \quad \text{as } n \to \infty \ . \]
\end{theorem}

We point out that the same result holds in the critical case (see \cite{agkv}), whereas it is no longer true in the moderate and strongly subcritical case (see e.g. \cite{gkv}). As a corollary we obtain from Proposition \ref{pr1} below the following result.

\begin{corollary}
Under A1 to A3 there is a number $0 < \kappa' < \infty$ such that
\[ \mathbb{P} \{ Z_n > 0\} \ \sim \ \kappa' \frac{\gamma^n}{na_n} \ . \]
\end{corollary}

The next theorem gives convergence of the laws of $Z_n$, conditioned on survival.

\begin{theorem} \label{th2}
Under A1 to A3 the conditional laws $\mathcal L(Z_n\, | \, Z_n > 0)$, $n \ge 1$, converge weakly to some probability distribution on the natural numbers. Moreover the sequence $\mathbb{E} [Z_n^\vartheta \, | \, Z_n > 0]$ is bounded for any $\vartheta < \beta$, implying convergence to the corresponding moment of the limit distribution.
\end{theorem}

Our last theorem describes the limiting behavior of the rescaled generation size process $e^{-S_k}Z_k$ for $r_n \le k \le n-r_n$, where $(r_n)$ is a sequence of natural numbers with $r_n \to \infty$ (and certainly $r_n< n/2$). Thus we consider the process $Y^n=\{Y^n_t, t \in [0,1]\}$, given by
\[ Y^n_t := \exp(-S_{r_n+\lfloor (n-2r_n)t\rfloor} )Z_{r_n+\lfloor (n-2r_n)t\rfloor} \ . \]
This process has asymptotically paths of a constant random value. More precisely:
\begin{theorem} \label{th3}
Under A1 to A3, there is a process $\{W_t, t \in [0,1]\}$ such that as $n\rightarrow\infty$ 
\[
\mathcal L \big(Y_t^n, t \in [0,1]\ \big| \ Z_n>0\big)\ \Rightarrow \ \mathcal L \big(W_t , \ t \in [0,1] \big)
\]
weakly in the Skorohod space $D[0,1]$. Moreover, there is a random variable $W$ such that $W_t=W$ a.s. for all $t\in[0,1]$ and
\[
\mathbb{P}\{0<W<\infty\} \ = \ 1 \ . 
\]
\end{theorem}
\noindent
Weaker versions of this results can be found in \cite{af_98,gkv}.

Thus we have the following scenario in the weakly subcritical case (being different from other cases): Given $Z_n > 0$ the value of $Z_k$ is of bounded order for $k$ close to $0$ and close to $n$. Inbetween $S_k$ takes large values, as can be seen from the proofs. In the first part, roughly up to time $\lfloor \epsilon n\rfloor$, $S_k$ is increasing exponentially fast, and the growth of $Z_k$ resembles that of supercritical growth. Then $Z_k$ follows the value of $e^{S_k} = \mathbb E[Z_k \, | \, \Pi]$ in a completely deterministic manner, up to a random factor $W > 0$. Afterwards this behaviour persists as long as $S_k$ remains large. Only at the end $S_k$ returns to 0 in the manner of a random walk excursion (as in \cite{ig}), and $Z_k$ is no longer tied to $S_k$. For further explanations we refer to \cite{bgk}. 

For the proof we develop several limit theorems for random walks $(S_n)_{n\in\mathbb{N}_0}$, conditioned to stay positive up to time $n$, for functionals, which depend primarily on the values of $S_k$ with $k$ being close to $0$ or to $n$. These theorems are presented in the following section. The proofs of the theorems are given in the closing section.

\section{Some limit theorems for random walks}      
\setcounter{equation}{0} 
In this section, we develop conditional limit theorems for a class of oscillating random walks without refering to branching processes. \\

Let $X_1,X_2,\ldots$ be a sequence of i.i.d. real-valued random variables and $S_0$ independent of $X_1,X_2,\ldots$. The random walk $S=(S_n)_{n\in\mathbb{N}_0}$ is defined by
\[S_n=S_0+X_1\ +\cdots + \ X_n\ .\] Our results are valid under a more general condition than A2, namely:

\paragraph{Assumption B.} {\em There are numbers $a_n>0$ such that $S_n/a_n$ converges in distribution to a law which is neither concentrated on $\mathbb{R}^+$ nor on $\mathbb{R}^-$.} \\

As is well-known the limit distribution is strictly stable with index $\alpha\in(0,2]$ with a density $s(x)$, such that $s(0)>0$.

Our theorems rely on conclusions from the theory of random walks, which we put together in this section. They rest on and substantially extend results due to Afanasyev \cite{af_90}, Bertoin and Doney \cite{bedo}, Hirano \cite{hi}, Iglehart \cite {ig}, Keener \cite{ke}, and others. 

In the sequel we shall also consider the possibility that the random walk starts from any point $x \in \mathbb{R}$ or from an initial distribution $\mu$. In such cases we write for probabilities as usual $\mathbf{P}_x\{\cdot\}$ or $\mathbf{P}_\mu\{\cdot\}$. We write $\mathbf{P}$ instead of $\mathbf{P}_0$. 

Duality will be an important tool later. Recall that given $n$ one may consider the dual objects $Q_i' := Q_{n-i+1}$, $X_i' := X_{n-i+1}$ for $i=1,\ldots,n$. Then the dual random walk is given by $S_i' := X_1'+\cdots + X_i'= S_n-S_{n-i}$, $S'_0=0$. We refrain from indicating the dependence on $n$ in the notation.

Let us introduce 
\[  M_n \ := \ \max(S_1, \ldots,S_n) \ , \quad L_n \ := \ \min(S_1, \ldots, S_n)\]
and the right-continuous functions $u: \mathbb{R} \to \mathbb{R}$ and $v: \mathbb{R} \to \mathbb{R}$ given by
\begin{align*} u(x)  \ &:= \ 1 + \sum_{k=1}^\infty \mathbf{P} \{-S_k\le x, M_k < 0 \} \ , \quad x \ge 0 \ ,\\
v(x) \ &:= \ 1 + \sum_{k=1}^\infty \mathbf{P} \{-S_k > x, L_k \ge 0\} \ , \quad x \le 0 \  
\end{align*}
and 0 elsewhere. 
In particular $u(0)=v(0)=1$. It is well-known that $u(x)=O(x)$, $v(-x)=O(x)$ for $x \to \infty$.

\paragraph{2.1 Large deviations for random walks.} The following precise large deviation estimates are extensions of known results. Some related results can be found in \cite{do2010}. Recall that $s(x)$ denotes the limiting density of $S_n/a_n$ and that  $s(0)>0$ under Assumption B. 

\begin{proposition} \label{pr1} For $\theta > 0,x \ge 0,$
\[\mathbf{E}_x [ e^{- \theta S_n} \, ;  \, L_n \geq 0] \ \sim \  s(0) b_n u(x) \int_0^\infty e^{-\theta z} v(-z) \, dz\ , \]
and for $\theta>0, x \le 0$ 
\[ \mathbf{E}_x[ e^{ \theta S_n} \, ;  \, M_n < 0] \ \sim \ s(0) b_n v(x)   \int^{\infty}_0 e^{-\theta z} u(z) \, dz\ , \]
with \[ b_n= (a_n n)^{-1} \ .\]
\end{proposition}

\noindent
For the proof we need the following lemma.

\begin{lemma}\label{le1} Let $(\beta_n)$ be a  regularly varying sequence with $\sum_{k=0}^{\infty} \beta_k < \infty$ and $d,e>0$. 
\begin{enumerate}[i)]
\item If $\delta_n \sim d \beta_n$, $\eta_n \sim e \beta_n$, then $\sum_{i=0}^n \delta_i \eta_{n-i} \sim c \beta_n$ with $c := d\sum_{k=0}^\infty \eta_k + e \sum_{k=0}^\infty \delta_k$ as $n\rightarrow\infty$.
\item  If $\sum_{k=0}^\infty \alpha_k t^k = \exp \big( \sum_{k=0}^\infty \beta_k t^k \big)$ for $|t|<1$, then $\alpha_n \sim  c \beta_n$ with $c :=   \sum_{k=0}^\infty \alpha_k$ as $n\rightarrow\infty$. 
\end{enumerate}
\end{lemma}

\noindent
{\em Proof.} i) is a well-known elementary fact and ii) is a special case of Theorem 1 in \cite{cnw}. \hfill $\Box$\\

\noindent
{\em Proof of Proposition \ref{pr1}.} Both claims are proven along the same lines. Since the first one has been considered (under stronger conditions) by Hirano \cite{hi}, let us turn to the second statement. By Assumption B and Stone's Local Limit Theorem (cf. \cite{bi}, section 8.4.1) for any interval $I$ of length $|I|$
\[ a_n \mathbf{P} \{ S_n \in I \} \ \to \ s(0) |I| \ .\]
Moreover the local limit theorem implies that there is a $c > 0$ such that
\begin{equation}
\mathbf{P} \{ S_n \in I \} \ \leq \ c/a_n  \label{local}
\end{equation}
uniformly in $n$ and all intervals $I$ of length at most $1$. Therefore for $\theta > 0$
\begin{align}
\sup_n a_n\mathbf{E}[ e^{\theta S_n}; S_n < 0] \ \le \ \sup_n \sum_{k = 0} ^\infty a_n\mathbf{P}\{ -k-1 \le  S_n < -k\} e^{-\theta k} \ < \ \infty \ . \label{2005}
\end{align}
Also for any $h>0$, 
\begin{align*}
a_n \sum_{k=0}^{\infty} e^{-\theta (k+1) h} \mathbf{P}&\{-(k+1)h\leq S_n <-kh\} \leq a_n\mathbf{E}[ e^{\theta S_n}; S_n < 0]
\\
&\leq a_n \sum_{k=0}^{\infty} e^{-\theta k h} \mathbf{P}\{-(k+1)h\leq S_n <-kh\}  \ .
\end{align*}
Now taking the limit $n\rightarrow \infty$, the limit and the sums interchange due to (\ref{2005}) and dominated convergence. Then taking the limit $h\rightarrow 0$ yields
\[ a_n \mathbf{E} [ e^{\theta S_n}; S_n <0] \ \to \ s(0)\int^{\infty}_0 e^{-\theta z} dz \ = \ \frac{s(0)} \theta \ .\]
Next the Baxter identity says that, for $|t|<1$ and $\theta>0$
\[ 1 + \sum_{k=1}^\infty t^k \, \mathbf{E} [e^{\theta S_k}; M_k < 0] \ = \ \exp \Big( \sum_{k=1}^\infty \frac{t^k}{k} \mathbf{E} [e^{\theta S_k}; S_k < 0] \Big)  \]
(cf. \cite{fe}, chapter XVIII.3 or \cite{bi}, chapter 8.9). Also $ \sum_{k=1}^\infty k^{-1} \, \mathbf{E} [e^{\theta S_k}; S_k < 0] < \infty$. 
Thus from Lemma \ref{le1} ii) it follows
\begin{align} na_n \mathbf{E} [e^{\theta S_n}; M_n < 0] \ &\to \ \frac {s(0)}\theta \Big( 1 + \sum_{k=1}^\infty \mathbf{E} [e^{\theta S_k}; M_k< 0]\Big) \notag \\
&= \ s(0)\int_0^\infty e^{-\theta z} u(z) \, dz  \ . \label{Hirano}
\end{align} 
In much the same way Hirano \cite{hi} obtained
\[na_n \mathbf{E} [e^{-\theta S_n}; L_n \geq 0] \ \to\ s(0)\int_0^\infty e^{-\theta z} v(-z) \, dz  \ .\]
by means of the corresponding Baxter identity ($L_n \ge 0$ and \mbox{$S_n \ge 0$} replacing $M_n < 0$ and $S_n < 0$). 
By the continuity theorem for Laplace transforms this generalizes to
\begin{align} 
na_n \mathbf{E} [e^{\theta S_n}; M_n < 0\, , \, S_n > -x] \ &\to \ s(0) \int_0^x e^{-\theta z} u(z) \, dz \ , \label{Hirano4}\\
na_n \mathbf{E} [e^{-\theta S_n}; L_n \ge 0\, , \, S_n < x] \ &\to \ s(0) \int_0^x e^{-\theta z} v(-z) \, dz \ , \label{Hirano2}
\end{align}
which for finite $x \geq 0$ now holds for every $\theta \in \mathbb{R}$. Note that the limit measures involved here have densities with respect to the Lebesgue measure and thus have no point masses, so that the convergence holds for any finite $x\geq 0$. 

Next let $x < 0$. By means of duality
\begin{align*} \mathbf{E} [e^{\theta S_n}; M_n < -x]  \ &= \ \sum_{i=0}^{n-1} \mathbf{E} [e^{\theta S_n}; S_0,\ldots,S_{i} \le S_i <-x \, , \, S_i > S_{i+1}, \ldots,S_n] \\
&\quad \quad + \ \mathbf{E}[e^{\theta S_n}; S_0,\ldots,S_n\leq S_n<-x] \\
&= \ \sum_{i=0}^n \mathbf{E} [e^{\theta S_i} ; L_i \ge 0, S_i <- x] \cdot \mathbf{E}[ e^{ \theta S_{n-i}}; M_{n-i} < 0] \ .
\end{align*}
This formula  together with (\ref{Hirano}), (\ref{Hirano2}) and the equations (note that $v(-z)$ is left continuous for $z>0$ and that $v(0)=v(0-)=1$)
\begin{align*}
1 + \sum_{k=1}^\infty \mathbf{E}[& e^{\theta S_k} ; L_k \ge 0,S_k < -x]  \\
&=\ 1+  \int_{(0,-x)} e^{\theta z} \, dv(-z) \ = \  e^{-\theta x}v(x)-\theta\int_0^{-x} e^{\theta z} v(-z) \, dz \ , \\
1+ \sum_{k=1}^\infty \mathbf{E}[&e^{\theta S_k} ; M_k < 0] \ = \ \theta  \int_0^\infty e^{-\theta z} u(z) \, dz 
\end{align*}
imply by means of Lemma \ref{le1} i) for $\theta > 0$ and $x<0$
\begin{equation}
na_n \mathbf{E} [e^{\theta S_n}; M_n < -x]  \ \to \ s(0)e^{-\theta x}v(x) \int_0^\infty e^{-\theta z} u(z) \, dz  \ ,\nonumber 
\end{equation}
which is equivalent to our claim. \hfill $\Box$\\

Related to these results are the following upper estimates. Their proofs shed some light on how the term $b_n=(a_nn)^{-1}$ comes into play. 

\begin{proposition} There is a number $c>0$ such that uniformly for all $x,y\ge 0$ and all $n$
\[ \mathbf{P}_x \{ L_n \ge 0, y-1 \le S_n <y \} \ \leq \ c \, b_n \,  u(x)v(-y) \ , \]
whereas for $x,y \le 0$
\[ \mathbf{P}_x \{ M_n < 0, y\le  S_n < y+1 \} \ \leq \ c\,  b_n \, v(x)u(-y) \ . \]
\end{proposition}

\noindent
{\em Proof.} As before, we prove the latter statement. Let $S'$ be the dual random walk and $L_i'$, $i=1,\ldots,n$,  the corresponding minima. Denote
\begin{align*} A_n \ := \ \{&M_{\lfloor n/3 \rfloor}< -x\} \\ A_n' \ := \ \{&L_{\lfloor n/3 \rfloor}' \ge y\} \ , \\
A_n'' \ := \ \{&y-x \le S_n < y-x+1\} \\ 
= \ \{&y-x  -T_n \le S_{\lfloor 2n/3 \rfloor} - S_{\lfloor n/3 \rfloor } < y-x-T_n+1 \} \ , 
\end{align*}
with $T_n := S_{\lfloor n/3 \rfloor} + S_n - S_{\lfloor 2n/3 \rfloor}$. Let $\mathcal{A}_n$ be the $\sigma$--field generated by $X_1,\ldots,X_{\lfloor n/3 \rfloor}$ and $X_{\lfloor 2n/3 \rfloor +1}, \ldots,X_n$. Then $T_n$ is $\mathcal{A}_n$--measurable, whereas $S_{\lfloor 2n/3 \rfloor} - S_{\lfloor n/3 \rfloor }$ is independent of $\mathcal{A}_n$, consequently from (\ref{local}) and the fact that $(a_n)$ is regularly varying there is a $c>0$ such that
\[ \mathbf{P} \{ A_n'' \, | \, \mathcal{A}_n \} \ \le \ c a_n^{-1} \ . \]
Since $A_n,A_n'$ are $\mathcal{A}_n$-measurable and independent, it follows
\[ \mathbf{P} \{ A_n \cap A_n' \cap A_n'' \} \ \le \ c a_n^{-1} \mathbf{P}\{A_n\} \mathbf{P}\{A_n' \} \ . \]
Moreover from Lemma 2.1 in \cite{agkv} there is a number $\rho \in (0,1)$ and a slowly varying sequence $l_1,l_2,\ldots$ such that
\[ \mathbf{P}\{ L_n \ge y \} \ \le \ c_1 u(-y) n^{-\rho} l_n \ , \quad  \mathbf{P}\{ M_n < -x \} \ \le \ c_2 v(x) n^{\rho -1}l_n^{-1} \ , \]
and we end up with the uniform estimate
\[ \mathbf{P} \{ A_n \cap A_n' \cap A_n'' \}  \ \le \  c v(x)u(-y) \, b_n  \]
for $c$ sufficiently large. Now since $M_{\lfloor n/3\rfloor}\leq M_n$ and $L'_{\lfloor n/3\rfloor} \leq S_n-M_n$,
\[ \{ M_n < -x, y-x \le S_n < y-x +1 \} \ \subset \ A_n \cap A_n' \cap A_n''  \ , \]
and the claim follows. \hfill $\Box$ 

\begin{corollary} \label{Korollar}
For any $\theta>0$ there is a $c>0$ (depending on $\theta$) such that for all $x,y\geq 0$
\[
\mathbf{E}_x\big[e^{-\theta S_n}; L_n\geq 0, S_n\geq y\big ] \leq c \ b_n  u(x)  v(-y) \ e^{-\theta y}
\]
and for all $x,y\leq 0$
\[
\mathbf{E}_x\big[e^{\theta S_n}; M_n< 0, S_n< y\big ] \leq c \ b_n  u(-y)  v(x) \ e^{\theta y} \ .
\]\end{corollary}
\noindent
{\em Proof.} Again we consider the latter statement. Let $\theta>0$ and $x,y\leq 0$. We use the inequalities
\begin{align*}
u(x+y)&\leq u(x)+u(y) \ , \\
u(x+y)&\leq 2u(x)u(y) \ .
\end{align*}
The first inequality is a consequence of the representation of $u$ as renewal function ($u(x)$ is the expected number of ladder points in the interval $[0,-x]$ plus one; see \cite{fe}, chapter XII). The second inequality follows directly from the first one (as $u(x)\geq 1$ for all $x\geq 0$).
Then
\begin{align*}
\mathbf{E}_x\big [ e^{\theta S_n};& \ M_n<0,S_n<y\big ] \nonumber \\
	&\leq \ \sum_{k=1}^{\infty} e^{\theta(y-k)} \mathbf{P}_x\big\{ M_n<0,-k+y\leq S_n < -k+y+1\big\} \\
	&\leq \ e^{\theta y} \sum_{k=1}^{\infty} e^{-\theta k} c \ b_n v(x) u(k-y-1) \\
	&\leq \ 2c \ b_n e^{\theta y} v(x) u(-y) \sum_{k=0}^{\infty} e^{-\theta (k+1)} u(k) \ . \tag*{\qed}
\end{align*} 

\paragraph{2.2 The probability measures $\mathbf{P}^+$ and $\mathbf{P}^-$.} The fundamental properties of $u,v$ are the identities
\begin{equation} \begin{array}{rl} \mathbf{E}  [u(x+X); X + x \ge 0] \ = \ u(x) \ ,  &x \ge 0 \ ,  \\   \mathbf{E} [ v(x+X);X+x<0] \ = \ v(x) \ ,   &x \le 0 \ , \end{array}  \label{harm}
\end{equation}
which hold for any oscillating random walk. 


We use them to introduce the probability measures $\mathbf{P}^+$ and $\mathbf{P}^-$. The construction procedure is standard and explained for $\mathbf{P}^+$ in detail in \cite{agkv} and \cite{bedo}. 

The probability measures $\mathbf{P}^+_x$, $x \ge 0$ are defined as follows. Assume that the random walk $(S_n)_{n\in\mathbb{N}_0}$ is adapted to some filtration $\mathcal F=(\mathcal F_n)$ and that $X_{n+1}$ is independent ot $\mathcal F_n$ for all $n \ge 0$. For every sequence $R_0,R_1, \ldots$ of $\mathcal S$-valued random variables, adapted to $\mathcal F$ and every integrable function $g:\mathcal{S}^{n+1}\rightarrow \mathbb{R}$, $n\in\mathbb{N}$, $\mathbf{P}^+_x$ fulfills
\begin{align*} \mathbf E^+_x[g(R_0,\ldots,R_n)] \ = \ \frac{1}{u(x)}\mathbf{E}_x[g(R_0,\ldots,R_n)u(S_n); L_n \ge 0] \ , \ n\in\mathbb{N}_0 \ . 
\end{align*}
This is the Doob transform from the theory of Markov chains. In particular, under $\mathbf{P}^+$ $S_0,S_1,\ldots$ is a Markov process with state space $[0,\infty)$ and transition probabilities 
\[ P^{\scriptscriptstyle{+}}(x,dy) \ := \ \frac {1}{u(x)} \mathbf{P}\{x+X \in dy\} u(y)1_{\{y\geq 0\}} \ , \quad x \ge 0 \ .
\]
It is the random walk conditioned never to enter $(-\infty,0)$.

Similarly $v$ gives rise to probability measures $\mathbf{P}^-_x$, $x \le 0 $, characterized by the equation
\begin{align*} \mathbf E^-_x[g(R_0,\ldots,R_n)] \ = \ \frac{1}{v(x)}\mathbf{E}_x[g(R_0,\ldots,R_n)v(S_n); M_n < 0] \ . 
\end{align*}
Under $\mathbf{P}^-_x$ the process $S_0,S_1, \ldots$ becomes a Markov chain with state space $\mathbb{R}^-$ and transition kernel
\[ P^{\scriptscriptstyle{-}}(x,dy) \ := \ \frac {1}{v(x)} \mathbf{P}\{x+X \in dy\} v(y)1_{\{y<0\}} \ , \quad x \le 0 \ .
\]
Note that $P^{\scriptscriptstyle{-}}(x, [0,\infty)) =0$, thus the Markov process never enters $[0,\infty)$ again. It may, however, start from the boundary $x=0$. Intuitively it is the random walk conditioned never to return to $[0,\infty)$.

\paragraph{Remark.} Under $\mathbf{P}^+_x$ the process $(S)_{n\in\mathbb{N}_0}$ may return to 0, however, under $\mathbf{P}^-_x$ this possibility is excluded. We remark that this subtlety has little impact: For $x<0$ there is a difference only for those $x$, where $v(x) \neq v(x-)$, that is for at most countably many $x$. In particular no difference occurs, if one considers (as in the sequel) measures $\mathbf{P}^-_\nu$ having an initial distribution $\nu$ without atoms.\hfill $\Box$ \\

\paragraph{2.3 Some conditional limit theorems.} By means of the measures $\mathbf{P}^+_x$, $\mathbf{P}^-_x$ we now generalize a result due to Hirano \cite{hi} on the limit behavior of certain conditional distributions. For $\theta>0$, let $\mu_{\theta}$, $\nu_{\theta}$ be the probability measures on $\mathbb{R}^+$ and $\mathbb{R}^-$ given by their densities
\[ \mu_{\theta}(dz) \ := \ c_1e^{ -\theta z} u(z)1_{\{z \ge 0\}} \, dz  \ , \quad \nu_{\theta}(dz) \ := \ c_2 e^{\theta z} v(z)1_{\{z < 0\}} \, dz  \]
with $c_1^{-1}=c_{1\theta}^{-1} = \int_0^\infty e^{ - \theta z} u(z) \, dz$, $c_2^{-1}=c_{2\theta}^{-1} = \int_{-\infty}^0 e^{\theta  z} v(z) \, dz$. 

As above let $R_0,R_1, \ldots$ be a sequence of $\mathcal S$-valued random variables, adapted to $\mathcal F$. Also let $Q_1,Q_2, \ldots$ be a sequence of i.i.d. random variables with values in some space $\mathcal D$ and adapted to $\mathcal F$, such that $Q_{n+1}$ is independent of $\mathcal F_n$ for all $n \ge 0$. Additionally, assume that $X_i$ is $\sigma(Q_i)$-measurable for all $i\geq 1$.

\begin{proposition} \label{pr2} For given $\theta>0$ and $i,j \ge 0$ let $U:=g(R_0,\ldots,R_i)$ and $V:=h(Q_1,\ldots,Q_j)$ be real-valued, bounded random variables with suitable bounded, measurable functions $g: \mathcal{S}^{i+1} \to \mathbb R$, $h:\mathcal{D}^{j} \to \mathbb R$. Also let $\varphi: \mathbb R \to \mathbb R$ be bounded and continuous. Denote $\tilde V_n:=h(Q_{n},\ldots,Q_{n-j+1})$. Then for $x \geq 0$ 
\begin{align*}
\frac{\mathbf{E}_x\big[U\tilde V_n\varphi(S_{n}) e^{- \theta S_n}\, ; \, L_n \ge 0\big]}{\mathbf{E}_x[e^{-\theta S_n} \, ; \, L_n \ge 0]} \ \to \ \mathbf E^+_x [U] \mathbf E^-_{\nu_{\theta}} \big[V\varphi(-S_0)\big]  \end{align*}
and for $x \le 0$
\begin{align*}
\frac{\mathbf{E}_x[U \tilde V_n\varphi(S_{n})e^{\theta S_n} \, ; \, M_n < 0]}{ \mathbf{E}_x[e^{\theta S_n} \, ; \, M_n < 0]} \ \to \ \mathbf E^-_{x} [U] \mathbf E^+_{\mu_{\theta}} [V\varphi(-S_0)]  \ .
\end{align*}
\end{proposition}

\noindent
{\em Proof.} The proofs of both claims are similar. From Proposition \ref{pr1} for $\lambda \ge 0$, $x,y \le 0$
\[ \frac{\mathbf{E}_{y}[e^{ (\lambda+\theta) S_{n-j}} ; M_{n-j}<0] } {\mathbf{E}_{x}[e^{\theta S_n} ; M_n<0]}
 \ \to \ \frac{v(y)\int_0^\infty e^{-(\lambda+\theta)z} u(z) \, dz}{ v(x)\int_0^\infty e^{-\theta z} u(z) \, dz} \ ,
\]
consequently by the continuity theorem for Laplace-transforms for $\varphi :\mathbb{R} \to \mathbb{R}$, bounded and a.s. continuous with respect to $\mu_{\theta}$
\[\frac{\mathbf{E}_{y}[\varphi(S_{n-j})e^{ \theta S_{n-j}} ; M_{n-j}<0] } {\mathbf{E}_{x}[e^{ \theta S_n} ; M_n<0]}
 \ \to \ \frac{v(y)}{ v(x)} \int \varphi(-z) \, \mu_{\theta}(dz) \ . \]
In particular this proves the proposition for $i=j=0$. Note that if $\varphi :\mathbb{R} \to \mathbb{R}$ is positive and a.s. continuous but possibly no longer bounded, we may conclude by a truncation procedure that
\begin{equation} \liminf_n \frac{\mathbf{E}_{y}[\varphi (S_{n-j})e^{ \theta S_{n-j}} ; M_{n-j}<0] } {\mathbf{E}_{x}[e^{ \theta S_n} ; M_n<0]}
 \ \ge \ \frac{v(y)}{ v(x)} \int \varphi (-z) \,  \mu_{\theta}(dz) \ . \label{limitdistr}
\end{equation}

In the general case let us assume without loss of generality $0 \le g,h \le 1$. From the Markov property for $n \ge i+j$
\begin{align*}
\mathbf{E}_{x}[U\tilde V_n\varphi(S_{n})e^{\theta S_{n}} \, ; \, M_n <  0] = \ \mathbf{E}_{x} \Big[U\psi_{n-i}(S_i) \, ; \, M_i < 0 \Big] \ ,
\end{align*}
where for $n\ge j$,
\[ \psi_n(y) \ := \ \mathbf{E}_{y} [\tilde V_n\varphi(S_{n})e^{\theta S_n} \, ; \, M_n < 0] \ . \]
By assumption, $\varphi$ is a bounded, continuous function. Thus,	 discontinuities of $\psi_j(y) =\mathbf{E}[\tilde V_j\varphi(S_j+y)e^{\theta (S_j+y)} \, ; \, M_j < -y]$ can only arise from discontinuities of $e: y\rightarrow \mathbf{P}\{M_j<y\}$. As bounded, monotone function, $e(\cdot)$ has at most countably many points of discontinuity. Thus the same holds for $\psi_j$ and $\psi_j$ is a.s. continuous with respect to $\mu_{\theta}$. Therefore it follows from 
\[\psi_{n-i}(y) \ = \ \mathbf{E}_y [ \psi_j(S_{n-i-j}) \, ; \, M_{n-i-j}< 0] 
\]
and from (\ref{limitdistr})
\[\liminf_n \frac{\psi_{n-i}(y) }{\mathbf{E}_{x} [e^{ \theta S_n} ; M_n<0]}\ \ge \ \frac{v(y)}{ v(x)} \int \psi_j(-z) \, e^{\theta z} \, \mu_{\theta}(dz) \ .
\]
By means of Fatou's Lemma 
\begin{align*}
&\liminf_n \mathbf{E}_x[U\tilde V_n\varphi(S_{n})e^{\theta S_n} ;  M_n < 0] \, \big/ \,  \mathbf{E}_x[e^{ \theta S_n} ; M_n<0] \\ & \qquad  \ge \ v(x)^{-1}\mathbf{E}_x[Uv(S_i);M_i<0] \cdot \int \psi_j(-z) \, e^{\theta z} \, \mu_{\theta}(dz)  \ . 
\end{align*} 
The first part of the righthand side is equal to $\mathbf E^-_x [U]$. As to the other part we use the duality transformation $Q'_i:= Q_{j-i+1}$, $i=1,\ldots,j$ and the corresponding path $S_1', \ldots,S_j'$, the invariance of the Lebesgue measure under the shift transformation $z \mapsto z+S_j'$ and the fact, that  the set $\{z : \min(S_0',\ldots,S_{j-1}') =-z\}$ has Lebesgue measure 0 and that $S_0'=0$, to obtain
\begin{align*}
&\int_0^\infty \psi_j(-z)u(z) \, dz  \\ 
&\ = \ \int_{0}^\infty \mathbf{E} [h(Q_1',\ldots,Q_j')\varphi(S_j'-z)e^{\theta(S_j'-z)}; S_0',\ldots,S_{j-1}' >S_j'-z ]u(z) \, dz \\
&\ = \ \mathbf{E} \int_{-S_j'}^\infty h(Q_1', \ldots, Q_j')\varphi(-z)e^{-\theta z} 1_{\{S_0'\ge -z\}}1_{\{S_1',\ldots,S_{j-1}' \ge -z\}}u(S_j'+z) \, dz \\
&\ = \ \mathbf{E} \int_{-S_0'}^\infty h(Q_1', \ldots, Q_j')\varphi(-z)e^{-\theta z} 1_{\{S_j'\ge -z\}}1_{\{S_1',\ldots,S_{j-1}' \ge -z\}}u(S_j'+z) \, dz \\
&\ = \ \mathbf{E} \int_{0}^\infty h(Q_1', \ldots, Q_j')\varphi(-z)e^{-\theta z} 1_{\{L_j' \ge -z\}} u(S_j'+z) \, dz \\
&\ = \ \int_{0}^\infty \mathbf{E}_{z} [h(Q_1,\ldots,Q_j)\varphi(-z)u(S_j); L_j \ge 0] e^{-\theta z} \, dz \\
&\ = \ \int_0^\infty \mathbf E^+_z [V\varphi(-S_0)] u(z) e^{-\theta z}\, dz\ .
\end{align*}
Altogether we end up with the estimate
\begin{align*}  
\liminf_n  \frac{\mathbf{E}_{x}[U\tilde V_n \varphi(S_{n})e^{\theta S_{n}} \, ; \, M_n <  0] }{ \mathbf{E}_{x}[e^{ \theta S_n} ; M_n<0]} \ge \ \mathbf E^-_x [U] \mathbf E^+_{\mu_{\theta}} [V\varphi(-S_0)] \ . \end{align*}
Finally replace in this estimate first $g$ by $1-g$ and $h$ by $1$ (i.e. $U$ by $1-U$ and $\tilde V_n,V$ by $1$) and second $h$ by $1-h$. Then these estimates altogether entail our claim. \hfill $\Box$ \\

We shall also use a dual version of the last proposition. Let
\begin{align}
 \tau_n \ := \ \min \{ i \le n : S_i  = \min(S_0, \ldots , S_n)\} \label{taun} 
\end{align}
be the moment of the first random walk minimum up to time $n$.

\begin{proposition} \label{pr21} Under the assumptions of  Proposition \ref{pr2}
\begin{align*}
\frac{\mathbf{E}[U\tilde V_n \varphi(S_n) e^{\theta S_n} \, ; \, \tau_n=n] }{ \mathbf{E}[e^{\theta S_n} \, ; \, \tau_n=n]}\ \to \ \mathbf E^+_{\mu_{\theta}} [U\varphi(-S_0)] \mathbf E^- [V]  \ .
\end{align*}
\end{proposition}

\noindent
{\em Proof.} There is a bounded, measurable function $\psi : \mathcal{D}^{i} \to \mathbb{R}$ such that a.s.
\[ \psi(Q_1,\ldots,Q_i) \ = \ \mathbf{E}[U |  Q_1,Q_2,\ldots, Q_i] \ .\]
By duality for $i+j \le n$
\begin{align*}
\mathbf{E}[U&\tilde V_n \varphi(S_n) e^{\theta S_n} \, ; \, \tau_n=n] \\ &= \ 
\mathbf{E}[\psi(Q_1,\ldots,Q_i)h(Q_{n}\ldots,Q_{n-j+1}) \varphi(S_n)e^{\theta S_n} \, ; \, \tau_n=n] \\
&= \ \mathbf{E}[h(Q_1',\ldots,Q_j')\psi(Q_n',\ldots,Q_{n-i+1}' )\varphi(S_n') e^{\theta S_n'} \; ; \; M_n'< 0] \ .
\end{align*}
Moreover
\[ \mathbf E^+_{\mu_{\theta}} [\psi(Q_1,\ldots,Q_i)\varphi(-S_0)] \ = \ \mathbf E^+_{\mu_{\theta}} [U\varphi(-S_0 )] \ ,\]
thus the claim follows from the preceding proposition.\qed\\

The next results on weak convergence generalizing the last propositions are in the spirit of Lemma 2.5 in \cite{agkv}.

\begin{theorem} \label{pr4}
Let $0 < \delta < 1$. Let $U_n=g_n(R_0,\ldots,R_{\lfloor \delta n\rfloor})$, $n\ge 1$, be  random variables with values in an Euclidean (or polish) space $S$ such that 
\[ U_n \ \to \ U_\infty   \quad \mathbf{P}^+ \text{-a.s.}  \]
for some $S$-valued random variable $U_\infty$. Also let $V_n=h_n (Q_1, \ldots, Q_{\lfloor \delta n\rfloor})$, $n\ge 1$, be random variables with values in an Euklidean (or polish) space $S'$ such that 
\[ V_n \ \to \ V_\infty   \quad \mathbf{P}^-_x \text{-a.s.}  \]
for all $x \le 0$ and some $S'$-valued random variable $V_\infty$. Denote
\[ \tilde{V}_n \ := \ h_{n}(Q_n, \ldots, Q_{n-\lfloor \delta n \rfloor+1}) \ . \]
Then for $\theta > 0$ and for any bounded, continuous function $\varphi:S\times S' \times \mathbb R \to \mathbb{R}$ as $n \rightarrow \infty$
\begin{align*} \mathbf{E} [ \varphi(U_n,&\tilde{V}_n,S_n) e^{-\theta S_n}\; ; \; L_n \ge 0] \; \big/ \; \mathbf{E}[e^{-\theta S_n};L_n \ge 0] \  \\ &\to \ \iiint \varphi(u,v,-z) \, \mathbf{P}^+\{U_\infty \in du\} \mathbf{P}^-_z\{V_\infty \in dv\} \nu_{\theta}(dz) \ . \end{align*} 
\end{theorem}

%

The following theorem is a counterpart.

\begin{theorem} \label{pr5}
Let $U_n, V_n,\tilde V_n$, $n=1,2,\ldots, \infty$ be as in be as in Theorem \ref{pr4}, now fulfilling
\[ U_n \ \to \ U_\infty   \quad \mathbf{P}^+_x \text{-a.s.} \ , \quad V_n \ \to \ V_\infty   \quad \mathbf{P}^- \text{-a.s.}\]
for all $x \ge 0$. Then for any bounded, continuous function $\varphi: S\times S' \times \mathbb R \to \mathbb{R}$ and for $\theta > 0$ as $n \rightarrow \infty$
\begin{align*} \mathbf{E} [ \varphi(U_n,&\tilde{V}_n,S_n)e^{\theta S_n} \; ; \;\tau_n=n] \; \big/\;  \mathbf{E}[e^{\theta S_n};\tau_n=n]  \\ &\to \ \iiint \varphi(u,v,-z) \, \mathbf{P}^+_z\{U_\infty \in du\} \mathbf{P}^-\{V_\infty \in dv\} \mu_{\theta}(dz) \ . \end{align*} 
\end{theorem}

The proofs of all three theorems are much the same. We prove the third one.\\

\noindent
{\em Proof of Theorem \ref{pr5}.} The proof relies on two estimates, which allow to switch from $\mathbf P$ to $\mathbf{P}^+$ resp. $\mathbf{P}^-$. First we look at the case $\varphi(u,v,z)= \varphi_1(u)$, where the function $\varphi_1$ is bounded by 1 and depends only on $u$. Then by the Markov property
\begin{align*} \mathbf{E}[\varphi_1(U_n)e^{\theta S_n}; \tau_n=n]  \ = \ \mathbf{E} [ \varphi_1(U_n) \psi_{n-\lfloor \delta n \rfloor  }(S_{\lfloor \delta n \rfloor }, \min(S_0,\ldots,S_{\lfloor \delta n \rfloor }))]
\end{align*}
with
\[ \psi_n(x,y) \ := \ \mathbf{E}_x[ e^{\theta S_n}; \tau_n=n, S_n < y] \ . \]
By duality and Corollary \ref{Korollar} for $x \ge y$
\begin{align*}
\psi_n(&x,y) \ = \ e^{\theta x} \mathbf{E} [e^{ \theta S_n}; \tau_n=n, S_n < y-x] \\
&= \ e^{\theta x}\mathbf{E}[ e^{\theta S_n}; M_n<0, S_n < y-x] \
\le \ c b_n u(x-y) e^{\theta y} \ ,
\end{align*}
and using $\min(S_0,\ldots,S_{\lfloor \delta n \rfloor })=L_{\lfloor \delta n \rfloor }\wedge 0$ it follows 
\begin{align*}
\big| \mathbf{E} [& \varphi_1({U}_n) e^{\theta S_n}; \tau_n =n] \big|\ \le \ \\
&  cb_{n-\lfloor \delta n \rfloor }  \mathbf{E}\big[| \varphi_1(U_n)| u\big(S_{\lfloor \delta n \rfloor }- L_{\lfloor \delta n \rfloor }\wedge 0\big) e^{\theta L_{\lfloor \delta n \rfloor}\wedge 0}\big] \ .
\end{align*}
By martingale property of $u$, (\ref{harm}), we have for any $y\in \mathbb{N}$ 
\begin{align*}
\mathbf{E}\big[u\big(&S_{\lfloor \delta n \rfloor }- L_{\lfloor \delta n \rfloor }\wedge 0\big) e^{\theta L_{\lfloor \delta n \rfloor }\wedge 0}; L_{\lfloor \delta n \rfloor } < -y\big] \\
&\le \ \sum_{j=y}^\infty e^{-\theta j} \mathbf{E} [u(S_{\lfloor \delta n \rfloor }+j+1); -j-1 \le L_{\lfloor \delta n \rfloor } < -j] \\
&\le \sum_{j=y}^\infty e^{-\theta j}  \mathbf{E}_{j+1} [u(S_{\lfloor \delta n \rfloor }) ; L_{\lfloor \delta n \rfloor } \ge 0]\\
&= \ \sum_{j=y}^\infty e^{-\theta j}u(j+1) \ .
\end{align*}
Also by duality and (\ref{Hirano})
\begin{align} \mathbf{E}[e^{\theta S_n}; \tau_n=n] \ = \ \mathbf{E}[e^{\theta S_n}; M_n<0] \ \sim \ c \ b_n \ .  \label{dual}
\end{align}
Thus, given $\varepsilon > 0$ and choosing $y$ sufficiently large, we have by Proposition \ref{pr1}
\begin{align*}
\big| \mathbf{E} [& \varphi_1({U}_n) e^{\theta S_n}; \tau_n =n] \big|\ \le \ \varepsilon \mathbf{E} [ e^{\theta S_n}; \tau_n =n] \\
&+  cb_{n-\lfloor \delta n \rfloor }  \mathbf{E}[| \varphi_1(U_n)| u(S_{\lfloor \delta n \rfloor }+y) ; L_{\lfloor \delta n \rfloor } \ge -y] \ .
\end{align*}
This leads to the estimate 
\begin{equation}
\frac{\big| \mathbf{E} [ \varphi_1({U}_n) e^{\theta S_n}; \tau_n=n] \big|}{\mathbf{E} [ e^{\theta S_n}; \tau_n =n]}\ \le \ \varepsilon  + c \mathbf E^+_{y}[| \varphi_1(U_n)|] \ , \label{estimate2}
\end{equation}
which for given $\varepsilon >0$ holds for $y$ sufficiently large.

Next we look at the case $\varphi(u,v,z)= \varphi_2(v)$, where $\varphi_2$ again is bounded  \linebreak by 1. By means of duality and the Markov property we obtain
\begin{align*}\mathbf{E} [\varphi_2(\tilde{V}_n)e^{\theta S_n} ; &\tau_n=n] \ = \  \mathbf{E} [\varphi_2(V_n)e^{\theta S_n}; M_n<0] \\
&= \ \mathbf{E} [\varphi_2(V_n) \psi_{n-\lfloor \delta n \rfloor }(S_{\lfloor \delta n \rfloor}); M_{\lfloor \delta n \rfloor} < 0] \ ,
\end{align*}
where
\[ \psi_n(x) \ := \ \mathbf{E}_x [e^{\theta S_n}; M_n < 0 ] \ . \]
By means of Corollary \ref{Korollar} there is a $c>0$ such that
\[ \big| \mathbf{E} [\varphi_2(V_n)e^{\theta S_n} ; M_n < 0]  \big| \
\le \ c b_{n-\lfloor \delta n \rfloor}  \mathbf{E} \Big[\big| \varphi_2(V_n)\big| v(S_{\lfloor \delta n \rfloor }) \, ; \, M_{\lfloor \delta n \rfloor } < 0 \Big] \ .\]
Recalling the definition of $\mathbf E^-$ and Proposition \ref{pr1} we obtain for a suitable $c>0$
\begin{equation}
\Big| \frac{\mathbf{E} [\varphi_2(\tilde{V}_{n})e^{\theta S_n} ; \tau_n =n] }{ \mathbf{E} [e^{\theta S_n} ; \tau_n=n] } \Big| \ \le \ c \mathbf E^- \big[ \big| \varphi_2(V_n)\big| \big] \ . \label{estimate1}\end{equation}

Now denote
\[\tilde{V}_{k,n} \ := \ h_k(Q_n,\ldots,Q_{n-{\lfloor \delta k\rfloor}+1})\ . \]
If $\varphi_3$ depends only on $z$ and is continuous and  bounded by 1, then  we obtain by means of (\ref{estimate2}) and (\ref{estimate1}) (replacing  $\varphi_1(U_n)$ and $\varphi_2(\tilde{V}_n)$ by $\varphi_1(U_n)-\varphi_1(U_k)$ and $\varphi_2(\tilde{V}_n)-\varphi_2(\tilde{V}_{k,n})$ in these estimates)
\begin{align*} &\frac{\big| \mathbf{E} [( \varphi_1(U_n) \varphi_2(\tilde{V}_n)-\varphi_1(U_k) \varphi_2(\tilde{V}_{k,n})) \varphi_3(S_n) e^{\theta S_n}; \tau_n=n]   \big|}{\mathbf{E} [ e^{\theta S_n}; \tau_n=n] } \\
&\quad \le  \ \frac{\mathbf{E} [| \varphi_1(U_n)-\varphi_1(U_k)|e^{\theta S_n}; \tau_n =n] }{\mathbf{E} [ e^{\theta S_n};\tau_n =n] }  + \frac{\mathbf{E} [| \varphi_2(\tilde{V}_n)-\varphi_2(\tilde{V}_{k,n})|e^{\theta S_n}; \tau_n=n] }{\mathbf{E} [ e^{\theta S_n}; \tau_n=n] }      \\
&\quad  \le  \ c \mathbf E^+_y[| \varphi_1(U_n)-\varphi_1(U_k)|] + \varepsilon + c\mathbf E^-[| \varphi_2(V_n)-\varphi_2(V_k)|] \ ,
\end{align*}
if $c,y$ are sufficiently large. Letting $n \to \infty$ we obtain by assumption and Proposition \ref{pr21}
\begin{align*}
\limsup_n \ &\Big|\frac{ \mathbf{E} [( \varphi_1(U_n) \varphi_2(\tilde{V}_n) \varphi_3(S_n) e^{\theta S_n}; \tau_n =n]}{\mathbf{E} [ e^{\theta S_n}; \tau_n=n] }\\ & \hspace{4cm}- \mathbf E^+_{\mu_{\theta}} [\varphi_1(U_k) \varphi_3(-S_0)]\mathbf E^- [\varphi_2(V_k)] \Big| \\
&  \le \ \varepsilon + c \mathbf E^+_y[| \varphi_1(U_\infty)-\varphi_1(U_k)|] + c\mathbf E^- [| \varphi_2(V_\infty)-\varphi_2(V_k)|] \ .
\end{align*}
Also by assumption the terms on the righthand side vanish for $k \rightarrow \infty$. Letting $\varepsilon \to 0$, our claim follows in the case $\varphi(u,v,z)= \varphi_1(u) \varphi_2(v) \varphi_3(z)$. As is well-known this case is sufficient for the proof of weak convergence. \qed

\section{Proof of theorems}
\setcounter{equation}{0} 

Define
\[ \eta_{i} \ := \ \sum_{y=0}^\infty y(y-1) Q_i(\{y\}) \Big/ \Big( \sum_{y=0}^\infty y Q_i(\{y\}) \Big)^2 \ , \quad i \geq 1 \ . \]

\begin{lemma} \label{le2} Assume A2 and  A3. Then for all $x \ge 0$
\[ \sum_{i=0}^\infty \eta_{i+1} e^{-S_i} \ < \ \infty \qquad \mathbf{P}^+_x \text{ -a.s.} \]
and for all $x \le 0$
\[ \sum_{i=1}^\infty \eta_{i} e^{S_i} \ < \ \infty \qquad \mathbf{P}^-_x \text{ -a.s.} \]
\end{lemma}

The proof of the first statement can be found in \cite{agkv} (see Lemma 2.7 therein under condition B1 and B2), the second one can be proven just the same way. \\

The branching mechanism can be neatly described by means of generating functions. Let
\[ f_j(s) \ := \ \sum_{i=0}^\infty s^i Q_j(\{i\}) \ , \quad 0 \le s \le 1 \ , \]
$j=1,2,\ldots$ and their compositions
\begin{align*}
f_{k,n} \ :=& \ f_{k+1}(f_{k+2}(\cdots f_n(s) \cdots)) \ , \quad 0 \le k < n \ , \\
f_{k,0} \ :=& \ f_k(f_{k-1}(\cdots f_1(s) \cdots)) \ , \quad 0 < k \ .
\end{align*}
As is well-known the branching property can be expressed as
\begin{equation} 
\mathbb{E} [s^{Z_n} \; | \;  \Pi,Z_k] \ = \ f_{k,n}(s)^{Z_k} \qquad \mathbb{P} \text{ -a.s.}
\end{equation}
and it also holds after a change of measure and conditioning (compare \linebreak section 3 in \cite{agkv}),
\begin{equation} 
\mathbf E^{\pm} [s^{Z_n} \; | \;  \Pi,Z_k] \ = \ f_{k,n}(s)^{Z_k} \qquad \mathbf P^{\pm} \text{ -a.s.}
\end{equation}

We shall use the following fact.

\begin{lemma} \label{lem}
For any $0 \le s \le 1$ the sequence $ f_{k,0}(s)^{\exp(-S_k)}$, $k\ge 1$, is non-decreasing and $ f_{k,0}(s)^{\exp(-S_k)} \ge s^{\exp(-S_0)}$.
\end{lemma}

\noindent
{\em Proof.} Without loss of generality assume $s >0$. We use the fact that the cumulant generating functions $c_k(\lambda):=\log f_k(e^{\lambda})$, $\lambda \le 0$, are convex. Since $c_k(0)=0$, this implies $c_k(\lambda) \ge c_k'(0)\lambda$ or, letting $\lambda =  \log t$, $0<t \le 1$,
\[ \log f_k(t) \ \ge\ f_k'(1) \log t = e^{X_k} \log t\ . \]
Choosing $t= f_{k-1,0}(s)$ and multiplying with $\exp(-S_k)$ gives the first statement. For $k=1$ the second inequality follows. \qed \\

Under the measure $\mathbb{P}$, Proposition \ref{pr1} translates to
\begin{corollary} \label{cor1}
For $x \ge 0,\theta > -\beta$
\[\mathbb{E}_x [ e^{- \theta S_n} \, ;  \, L_n \geq 0] \ \sim \  s(0) \gamma^n b_n u(x) e^{\beta x}\int_0^\infty e^{-(\theta + \beta)z} v(-z) \, dz\ , \]
and for $x \le 0$ and $\theta > \beta$ 
\[ \mathbb{E}_x[ e^{ \theta S_n} \, ;  \, M_n < 0] \ \sim \ s(0) \gamma^n b_n v(x)e^{\beta x}  \int^{\infty}_0 e^{-(\theta-\beta)z} u(z) \, dz\ . \]
\end{corollary}
\noindent
{\em Proof.} We only prove the second statement. Let $x\le 0$ and $\theta>\beta$. By the usual change of measure,
\begin{align*} 
  \mathbb{E}_x[e^{\theta S_n};& \ M_n<0] \  =  \  \gamma^n e^{\theta x} \mathbf{E}[e^{(\theta-\beta)S_n}; M_n<-x] \\
 = & \  \gamma^n e^{\beta x} \mathbf{E}_x[e^{(\theta-\beta)S_n}; M_n<0] \ 
\end{align*} 
and the result follows from Proposition \ref{pr1}. \qed \\

\noindent
With the corresponding change to $\theta=\vartheta+\beta$, in the following, Theorem \ref{pr4} and \ref{pr5} will be used. By the change of measure used before, the theorems can be applied to the measure $\mathbb{P}$. For later use we note that
\begin{align}
\mathbb E[e^{S_n}; \tau_n=n] \sim c \mathbb P\{L_n \ge 0\}  
\label{nochmal}
\end{align}
with some number $c >0$, which follows from the last corollary together with \eqref{dual}. 
\begin{lemma} \label{lem1}
Let $z \ge 1$ and let $m_n$, $n \ge 1$, be a sequence of natural numbers with $m_n \sim n/2$. Then, as $n \to \infty$, the random vector $(\exp(-S_{m_n})Z_{m_n},Z_n)$, given the event $\{Z_0=z, L_n \ge 0\}$,  converges in distribution to some random vector $(W,G)$ with values in $[0,\infty) \times \mathbb N_0$. Moreover the probability of  $G \ge 1$ is greater than $0$ and $W>0$ a.s. on the event $G \ge 1$.
\end{lemma}

\noindent
{\em Proof.} 
We prove convergence of $\mathbb{E} [U_n's^{Z_n} \, | \, Z_0=z, L_n \ge 0]$ for $0 < s \le 1$ and suitable bounded random variables $U_n'$.
Define
\[ \varphi(u',u'',v,x) \ := \   u'v^{u''\exp(x)} \ , \quad \text{ for }0 \le u' \le 1, u'' \ge 0, 0 \le v \le 1, x \in \mathbb R \ , \]
with $0^0=1$. For other values of $(u',u'',v,x)$ let  $ \varphi(u',u'',v,x)$ be such that $\varphi$ becomes a bounded, continuous function. In doing so points of discontinuity in $(u',0,0,x)$ are unavoidable, which will be bypassed in the sequel. Moreover let
\begin{align*} U_n=(U_n',U_n'') \ &:= \ (U_n',\exp(-S_{m_n}) Z_{m_n})  \ ,\\ \quad  V_n \ &:= \ f_{n-m_n,0}(s)^{\exp(-S_{n-m_n})} \ ,
\end{align*}
and thus $\tilde V_n = f_{m_n,n}(s)^{\exp(-(S_n-S_{m_n}))}$.
If we assume that $U_n'$ is a random variable with values between $0$ and $1$ of the form $U_n'= h(S_{m_n},Z_{m_n})$, then 
\[\mathbb{E}[U_n's^{Z_n} \, | \,  \Pi, Z_{m_n}] = U_n'f_{m_n,n}(s)^{Z_{m_n}} = \varphi(U_n, \tilde V_n, S_n)\] 
and
\[ \mathbb{E} [U_n's^{Z_n} \, | \, Z_0=z, L_n \ge 0] \ = \ \mathbb{E} [\varphi(U_n, \tilde V_n, S_n) \, | \, Z_0=z, L_n \ge 0] \ . \]

We would like to apply Theorem \ref{pr4}. From Lemma \ref{lem} it follows that $V_n$ converges to some random variable $V_\infty$ such that $0 < s^{\exp(-S_0)} \le V_\infty \le 1$.  From Proposition 3.1 in \cite{agkv} (with the Assumptions B1 and B2 therein) we see that  $U_n'' $  converges $\mathbf{P}^+$-a.s. to a random variable $U_\infty''$. Also $\mathbf{P}^+(U_\infty'' > 0) > 0$. Thus we just have to take care that $U_n'$  converges $\mathbf{P}^+$-a.s. to some random variable $U_\infty'$.

Now $\varphi $ is continuous in every point $(u',u'',v,x)$ with $v > 0$, thus we conclude from Theorem \ref{pr4} and standard  results on weak convergence that 
\begin{align}
 \mathbb{E} [U_n's^{Z_n} \, | \, Z_0=z, L_n \ge 0]\rightarrow \psi_z(s) \qquad, \ 0<s \le 1 \ ,\label{convhilf}
\end{align}
where
\begin{align*}  \psi_z(s)\ :=\ \iiint \varphi(u,v,-x) \, \mathbf{P}^+\{U_\infty \in du\} \mathbf{P}^-_x\{V_\infty \in dv\} \nu_{\beta}(dx) \ , \end{align*} 
and $u=(u',u'')$ and $U_\infty=(U_\infty',U_\infty'')$.
Note that the distribution of $U_\infty$ depends only on $z$ whereas the distribution of $V_\infty$ depends on $s$. From $s^{\exp(-S_0)} \le V_\infty \le 1$ it follows, that 
\[ \psi_z(s) \to \psi_z(1)=\mathbf E^+[U_\infty'] \ , \text{ as } s\to 1 \ . \]

First let us choose $U_n'=1$ for all $n$.  Then for $0 < s \le 1$ we obtain convergence of the generating function $\mathbb{E} [s^{Z_n} \, | \, Z_0=z, L_n \ge 0]$ to some function $\psi_z(s)$ with $\psi_z(s)\to 1$, as $s \to 1$. Thus $\mathcal L(Z_n\, | \, Z_0=z, L_n \ge 0)$ is weakly convergent to some probability measure with generating function $\psi_z$. In order to show that this measure is not the Dirac measure at $0$ we prove that $V_\infty < 1$ $\mathbf{P}^-_x$-a.s. for $s < 1$. To this end we use 
\[ - \log V_{n} \ge \exp(-S_{n-m_n})(1-f_{n-m_n,0}(s)) \]
together with an estimate for $f_{k,n}$ due to Agresti \cite{ag} (see also the proof of Proposition 3.1 in \cite{agkv}), which for $f_{k,0}$ reads
\[\exp(-S_k)(1-f_{k,0}(s)) \ge \Big( \frac{1}{1-s} + \sum_{i=1}^k \eta_{i} e^{S_i} \Big)^{-1} \ . \]
From Lemma \ref{le2} we see that $-\log V_\infty > 0$ and thus $V_\infty < 1$ $\mathbf{P}^-_x$-a.s.  Also, as already mentioned, $\mathbf{P}^+(U_\infty'' > 0) > 0$. By definition of $\varphi(u,v,x)$ this implies $\psi_z(s)< 1$ for $s < 1$. Therefore the corresponding probability measure is not concentrated at $0$.

Next we choose $s=1$ (thus $V_n=1$) and $U_n' := \chi(\exp(-S_{m_n})Z_{m_n})= \chi(U_n'')$, where $\chi: \mathbb R \to [0,1]$ denotes a continuous function. Then by (\ref{convhilf}), we obtain that
\[ \mathbb E\big[\chi(\exp(-S_{m_n})Z_{m_n})\, | \, Z_0=z, L_n\ge 0\big] \to \mathbf E^+[\chi(U_\infty'')] \ . \]
This gives weak convergence of $ \mathcal L\big(\exp(-S_{m_n})Z_{m_n}\, | \, Z_0=z, L_n\ge 0\big)$ to some probability measure on $\mathbb R^+$. The convergence in (\ref{convhilf}) also implies that
\[ \mathbb{E} [\chi(\exp(-S_{m_n})Z_{m_n})s^{Z_n} \, | \, Z_0=z, L_n \ge 0] \]
has a limit for any $0 < s \le 1$ and any bounded continuous $\chi$. Therefore, given the event $\{Z_0=z, L_n \ge 0\}$, the joint distribution of $\exp(-S_{m_n})Z_{m_n}$ and $Z_n$  is weakly convergent, too. We write the limiting distribution as the distribution of some pair $(W,G)$ of random variables with values in $\mathbb R^+ \times \mathbb N_0$ and, for ease of notation, we denote the corresponding probabilities and expectations by $\mathbf P$ and $\mathbf E$. We already proved that $\mathbf P\{ G \ge 1 \} > 0$.

For the last claim of the lemma we use our convergence result for $s=1$ and $U_n' :=  I_{\{Z_{m_n} \ge 1\}}\chi(\exp(-S_{m_n})Z_{m_n})$ with continuous $\chi$ with values in $[0,1]$. From Proposition 3.1 in \cite{agkv} it follows that $U_n' $ converges to $I_{\{U_\infty'' > 0\}} \chi(U_\infty'')$ $\mathbf P^+$-a.s. and consequently
\[ \mathbb E\big[I_{\{Z_{m_n} \ge 1\}}\chi(\exp(-S_{m_n})Z_{m_n})\, | \, Z_0=z, L_n\ge 0\big] \to \mathbf E [I_{\{W > 0\}} \chi(W)] \ . \]
On the other hand we know that
\[\mathbb E\big[I_{\{Z_{n} \ge 1\}}\chi(\exp(-S_{m_n})Z_{m_n})\, | \, Z_0=z, L_n\ge 0\big] \to \mathbf E [ \chi(W); G \ge 1] \ . \]
Now $I_{\{Z_n \ge 1\}} \le I_{\{Z_{m_n} \ge 1\}}$, therefore $\mathbf E [ \chi(W); G \ge 1] \le \mathbf E [I_{\{W > 0\}} \chi(W)]$. For $\eta>0$ (choosing an appropriate $\chi$) it follows that
\[ \mathbf P\{W=0, G \ge 1\} \le \mathbf P\{0 < W \le \eta\} \ . \]
Letting $\eta \to 0$, this gives $\mathbf P\{W=0, G \ge 1\}=0$, which is our last claim.
\qed

\begin{lemma} \label{lem2}
Let $m_n$, $n \ge 1$, be such that $m_n \sim n/2$ and $\tau_n$ be defined as in (\ref{taun}). Then the conditional distribution $\mathcal L\big((\exp(-S_{m_n})Z_{m_n},Z_n) \, | \, Z_n > 0, \tau_n=n\big)$ converges to some random vector $(W,G)$ with values in $(0,\infty)\times \mathbb N$. Moreover there is a number $0<\kappa <\infty$ such that
\[\mathbb{P} \{Z_n > 0, \tau_n=n \} \ \sim \ \kappa \mathbb{P}\{L_n \ge 0 \} \ . \]
\end{lemma}

\noindent
{\em Proof.} The proof is somewhat different from the preceding one.
For $a > 0$ let
\[ \varphi_a(u',u'',v,x) := u'(1-v^{u''\exp(x)})1_{\{ x \ge -a\}} e^{-x} \ , \]
for $0 \le u' \le 1$, $u'' \ge 0$, $0\le v \le 1$, $x \in \mathbb R$ and continue $\varphi_a$ to other values of $u',u'',v,x$ to a bounded, smooth function. Conditioning as above we obtain
\[\mathbb{E} \big[ U_n'(1- s^{Z_n}) 1_{\{S_n \ge -a\}}; \tau_n=n\big] \ =\ \mathbb{E} \big[ \varphi_a(U_n, \tilde V_n, S_n) e^{S_n} ; \tau_n=n] \ ,\]
where $U_n=(U_n',U_n'')$, $\tilde V_n$ are as in the last proof. Note, that the additional discontinuity at $x=-a$ has probability 0 with respect to the measure $\mu_{1-\beta}$. Thus we may apply Theorem \ref{pr5} to $\varphi_a(u,v,x)$. 

Moreover, as $1-s^{Z_n}\leq Z_n$ and by duality
\begin{align*}
\mathbb{E}\big[|U_n'|(1-s^{Z_n}) &1_{\{S_n < -a\}} ; \tau_n=n\big] \ \leq\ \mathbb{E}\big[Z_n;S_n < -a, \tau_n=n\big]  \\
& = \ \mathbb{E}\big[e^{S_n} ;S_n < -a , \tau_n=n\big] \ = \ \mathbb{E}\big[e^{S_n} ;S_n < -a , M_n < 0 \big] \ .
\end{align*}
In view of \eqref{Hirano4} (translated to the measure $\mathbb{P}$ by the usual tilting), there is for every $\varepsilon > 0$ an $a > 0$ such that
\[ \mathbb{E}\big[1-s^{Z_n} ;S_n < -a , \tau_n=n\big] \ \le\ \varepsilon \mathbb{E}\big[e^{S_n} ; \tau_n=n \big] \]
for all $n$. Hence we conclude that the statement of Theorem \ref{pr5} holds for $\varphi_\infty(u,v,x)$, too, and we obtain as in the preceding proof
\[ \hat \psi_n(s) \ :=\ \mathbb{E} \big[ U_n'(1- s^{Z_n}); \tau_n=n\big] \, / \, \mathbb{E}[e^{S_n}; \tau_n=n] \ \to \ \hat \psi(s) \]
for $0 < s \le1$, with
\[\hat \psi(s):= \iiint u'(1-v^{u''\exp(x)}) e^x  \, \mathbf{P}^+_x\{U_\infty \in du\} \mathbf{P}^-\{V_\infty \in dv\} \mu_{1-\beta}(dx)  \ . 
\]

First we note that $\hat \psi(s)$ is right continuous at $0$: With decreasing $s>0$ also the values of $V_\infty$ decrease, and  the integrand of $\hat \psi(s)$ increases. Also for $u',u'',x$ fixed the integrand is continuous in $s$ (with $0^0=1)$. Therefore monotone convergence implies $\hat \psi(0+)=\hat\psi(0)$.

Next we note that the functions $\hat\psi_n(s)$, $n \ge 1$, are uniformly bounded analytical functions on the complex unit disc and convergent for $0< s< 1$. As is well known this implies convergence of $\hat\psi_n(s)$ to an analytic function on the unit disc. In particular $\hat\psi_n(0)\to \hat \psi(0)$, since $\hat \psi(0+)=\hat\psi(0)$. Also this convergence implies that the coefficients of the power series $\hat\psi_n(s)$, namely
\[\frac{\mathbb{E}[U_n';Z_n>0,\tau_n=n]}{\mathbb{E}[e^{S_n};\tau_n=n]} \quad \text{and} \quad  \frac{\mathbb{E}[U_n';Z_n=k,\tau_n=n]}{\mathbb{E}[e^{S_n};\tau_n=n]} \text{ with } k \ge 1 \] 
are convergent for $n \to \infty$.

Now let us look at the case $U_n'=1$. Then we obtain the existence of the limits 
\[\kappa_0=\lim_{n\rightarrow\infty} \frac{\mathbb{P}\{Z_n>0,\tau_n=n\}}{\mathbb{E}[e^{S_n};\tau_n=n]} \quad \text{and} \quad  \lim_{n\rightarrow\infty} \frac{\mathbb{P}\{Z_n=k,\tau_n=n\}}{\mathbb{E}[e^{S_n};\tau_n=n]}  \] 
for $k \ge 1$.  
Also $\hat \psi(s) >0$ for $s < 1$, which follows exactly as $\psi_z(s)< 1$ in the proof of the last lemma.
This implies $\kappa_0>0$, which, together with \eqref{nochmal}, gives the last statement of the lemma. Also it follows that
\[\lim_{n \to \infty} \mathbb P\{Z_n=k \mid Z_n >0, \tau_n=n\} \]
exists for all $k \ge 1$. 
We have to verify that the limiting measure is a probability distribution, that is we have to prove that the sequence of conditional distributions of $Z_n$, given $\{ Z_n >0, \tau_n=n\}$, is tight. This follows from the estimate
\begin{align*} \mathbb{P} \{ Z_n > k, \tau_n=n\} \ &\le\ \frac 1k \mathbb{E} [Z_n; \tau_n=n] \ =\ \frac 1k \mathbb{E}[ e^{S_n} ; \tau_n=n] \\ &\sim \frac 1k \ \frac{\mathbb{P}\{Z_n>0;\tau_n=n\}}{\kappa_0}  \ .
\end{align*}
Thus all statements on $Z_n$ are proven.

Next we consider convergence of the conditional distribution $U_n''=\exp(-S_{m_n})Z_{m_n}$. For this purpose let $U_n'=I_{\{Z_{m_n}\ge 1\}}\chi(U_n'')$ with $\chi$ continuous and bounded with values between $0$ and $1$. From Proposition 3.1 in \cite{agkv} $U_n'$ converges to $I_{\{U_\infty''>0\}}\chi(U_\infty'')$ $\mathbf P_x^+$-a.s. for all $x$. Since $Z_n>0$ implies $Z_{m_n}>0$
\[ \frac{\mathbb E[\chi(U_n''); Z_n >0 , \tau_n=n]}{\mathbb{E}[e^{S_n};\tau_n=n]} = \frac{\mathbb E[U_n'; Z_n >0 , \tau_n=n]}{\mathbb{E}[e^{S_n};\tau_n=n]} = \hat \psi_n(0) \ . \]

From $\hat\psi_n(0) \to \hat \psi(0)$ and from the definition of $\kappa_0$ it follows
\begin{align*} \mathbb E[\chi(U_n''); \ &Z_n >0 , \tau_n=n] \, / \, \mathbb{P}\{Z_n>0,\tau_n=n\} \\ 
&\to 
\frac{1}{\kappa_0}\iiint1_{\{u''>0\}} \chi(u'')(1-v^{u''\exp(x)}) e^x  \,  \pi(du'',dv,dx) 
\end{align*}
with $\pi(du'',dv,dx)=\mathbf{P}^+_x\{U_\infty'' \in du''\}\mathbf{P}^-\{V_\infty \in dv\} \mu_{1-\beta}(dx)$. This implies weak convergence of the distribution of $U_n''$, given $\{Z_n >0, \tau_n=n\}$,  to a probability distribution. Also, because of the appearence of $1_{\{u''>0\}}$ in the integral this distribution is concentrated on $(0,\infty)$.
Finally we also have the convergence of the coefficients of the power series $\hat\psi_n(s)$ and consequently the existence  of the limits
\[ \lim_{n\rightarrow\infty} \frac{\mathbb{E}[\chi(U_n'');Z_n=k,\tau_n=n\}}{\mathbb{P}\{Z_n>0,\tau_n=n\} } \]
for $k \ge 1$. This implies convergence of the joint distribution of $\exp(-S_{m_n})Z_{m_n}$ and $Z_n$, given $\{Z_n >0, \tau_n=n\}$. 
\qed

\begin{lemma} \label{lem3}
Under A1 to A3, for every $B \subset \mathbb N=\{1,2,\ldots\}$ there exists $0\le \kappa(B) <\infty $ such that
\[ \frac{\mathbb{P} \{ Z_n \in B\}}{\mathbb{P} \{ L_n \ge 0 \}}\rightarrow \kappa(B)  \ . \]
Also $\kappa>0$ for $B = \mathbb N$.
\end{lemma}

\noindent
{\em Proof.}
We decompose at the moment, when $S_1,\ldots, S_n$ takes its minimum for the first time:
\[ \mathbb{P} \{Z_n\in B\} \ =\ \sum_{k=0}^n \mathbb{P}\big\{ Z_n\in B , \tau_k=k, \min_{k<l\le n} S_l \ge S_k\big\} \ . \]
Letting
\[ \xi_n(z) =  \mathbb{P}\{ Z_n\in B \, | \, Z_0=z, L_n \ge 0 \} \]
we obtain for fixed $m \ge 1$ and $n \ge 2m$
\begin{align} \mathbb{P} \{Z_n\in B\}\ &=\ \sum_{k=0}^{m-1} \mathbb{E} [\xi_{n-k}(Z_{k}); \tau_k=k] \mathbb{P}\{L_{n-k} \ge 0\} \label{deco1} \\
&\quad + \ \sum_{k=m}^{n-m} \mathbb{E} [\xi_{n-k}(Z_{k}); \tau_k=k] \mathbb{P}\{L_{n-k} \ge 0\} \label{deco2}\\
&\quad + \ \sum_{j=0}^{m-1} \mathbb{E} [\xi_{j}(Z_{n-j}); \tau_{n-j}=n-j] \mathbb{P}\{L_{j} \ge 0\} \label{deco3}\ .
\end{align}
As to the sum in \eqref{deco1}, $\xi_n(z)$ is bounded by $1$ and in view of Lemma \ref{lem1} converges for every $z \ge 1$. Also $\mathbb{P}\{L_{n-k} \ge 0\} \sim \gamma^{-k} \mathbb{P}\{L_n \ge 0\}$ by \linebreak Corollary \ref{cor1}. Therefore there is a number $0\le \kappa' <\infty$ (depending on $m$ and $B$) such that
\[ \sum_{k=0}^{m-1} \mathbb{E} [\xi_{n-k}(Z_{k}); \tau_k=k] \mathbb{P}\{L_{n-k} \ge 0\} \ = \ (\kappa'+o(1)) \mathbb{P}\{L_n \ge 0 \} \ . \]

The sum in \eqref{deco2} may be estimated from above by 
\[ \sum_{k=m}^{n-m} \mathbb{P}\{ Z_k > 0, \tau_k=k\} \mathbb{P} \{L_{n-k} \ge 0\} \]
or in view of Lemma \ref{lem2} by
\[ \sum_{k=m}^{n-m} \mathbb{P} \{L_{k}\ge 0\}  \mathbb{P} \{L_{n-k}\ge 0\} \]
up to some factor independent of $m$. In view of Corollary \ref{cor1} this may be bounded by
\[ \gamma^n \sum_{k=m}^{n-m} b_k\cdot b_{n-k} \ ,\]
again up to an uniform factor. Lemma \ref{le1} shows that this quantity is asymptotically equal to $\gamma^n b_n 2 \sum_{k=m}^\infty b_k$. Altogether, in view of Corollary \ref{cor1}  for any $\varepsilon > 0$ 
\[ \sum_{k=m}^{n-m} \mathbb{E} [\xi_{n-k}(Z_{n-k}); \tau_k=k] \mathbb{P}\{L_{n-k} \ge 0\} \ \le \ \varepsilon \mathbb{P} \{L_n \ge 0 \} \ , \]
if only $m$ is large enough.

Finally for the sum in \eqref{deco3} 
\begin{align*}&\mathbb{E} [\xi_{j}(Z_{n-j}); \tau_{n-j}=n-j] \\ &\qquad = \ \mathbb{E} [\xi_{j}(Z_{n-j})\, | \, Z_{n-j} >0,\tau_{n-j}=n-j] \mathbb{P}\{Z_{n-j} >0,\tau_{n-j}=n-j\} \ . 
\end{align*}
The first term in righthand-side of the above equation is bounded by 1 and it follows with Lemma \ref{lem2} that
\[ \sum_{j=0}^{m-1} \mathbb{E} [\xi_{j}(Z_{n-j}); \tau_{n-j}=n-j] \mathbb{P}\{L_{j} \ge 0\} \ = \ (\kappa''+o(1)) \mathbb{P} \{ L_n \ge 0\} \]
for some $\kappa'' \ge 0 $. 

Altogether, letting $m \to \infty$, these three statements imply the first claim of the lemma. Also, if $B = \mathbb N$, then because of Lemma \ref{lem1} the limit of $\xi_n(z)$ is strictly positive for all $z$. Because of Lemma \ref{lem2} $\kappa' > 0$ for $m$ sufficiently large. This gives the second statement. \qed

\paragraph{Proof of Theorem \ref{th1}.}  This theorem is contained in the last lemma.

\paragraph{Proof of Theorem \ref{th2}.} First we estimate $\mathbb{E} [Z_n^\vartheta ] $:
\begin{align*}
\mathbb{E} [Z_n^\vartheta ] \ &= \  \sum_{k=0}^n \mathbb{E}[ Z_n^\vartheta ; \tau_k=k, \min_{k<l\le n} S_l \ge S_k]\\
&= \ \sum_{k=0}^n \mathbb{E} [ \eta_{n-k}(Z_k)  ; \tau_k=k] \mathbb{P} \{ L_{n-k} \ge 0\} 
\end{align*}
with
\[ \eta_n(z) := \mathbb{E} [Z_n^\vartheta \, | \, Z_0=z, L_n \ge 0 ] \ . \]
From Jensen's inequality for $\vartheta<1$
\begin{align*}
\eta_n(z)\ &\le \ \mathbb{E} \big[ \mathbb{E}[Z_n\, | \,  \Pi]^\vartheta \, \big| Z_0=z, L_n \ge 0\big] \\
& = \ \mathbb{E} [ z^\vartheta \exp(\vartheta S_n) \, | \, L_n \ge 0] \ \le \ 
z\mathbb{E} [ \exp(\vartheta S_n) \, | \, L_n \ge 0] \ .
\end{align*}
In view of Corollary \ref{cor1} there is a $ c' > 0$ such that $
\eta_n(z) \ \le \ c'z $ for $\vartheta<\beta$
and $n \ge 1$, therefore
\[ \mathbb{E} [ \eta_{n-k}(Z_k)  ; \tau_k=k] \ \le \ c' \mathbb{E} [ Z_k ; \tau_k=k] \ = \ c' \mathbb{E} [ e^{S_k} ; \tau_k=k] \ \le \ c'' \mathbb{P}\{L_k \ge 0\} \ . \]
As in the proof of the last lemma this implies
\[ \mathbb{E} [Z_n^\vartheta ] \ \le \ c \mathbb{P} \{ L_n \ge 0 \} \ . \]
for a suitable $c > 0$. In view of Theorem \ref{th1} it follows that
$\mathbb{E} [Z_n^\vartheta \, | \, Z_n > 0]$ is bounded for $\vartheta < \beta$. 

This also gives tightness of the distributions of $Z_n$, given $Z_n > 0$. From the last lemma we see that 
$ \mathbb{P} \{ Z_n = a\, | \, Z_n > 0\} $
is convergent for $a \ge 1$. This completes the proof.
\qed
\paragraph{Proof of Theorem \ref{th3}}
First we consider 
\[
 Y^n_{1/2}= e^{-S_{\lfloor n/2\rfloor}} Z_{\lfloor n/2 \rfloor} \ . \]
We show that
\begin{align}\mathcal L\big( Y^n_{1/2}\ \big|\ Z_n>0\big) \stackrel{d}{\rightarrow} W \ , \label{0505}
\end{align}
where $W$ is an a.s. positive random variable. 

Let $\chi:\mathbb R \to \mathbb R$ be bounded and continuous.
As above in Lemma \ref{lem3} we consider the decomposition
\[
\mathbb{E}\big[\chi(Y^n_{1/2}); Z_n>0\big] = \sum_{k=0}^n \mathbb{E}\big[\chi(Y^n_{1/2}); Z_n>0 , \tau_k=k, \min_{k<l\le n} S_l \ge S_k\big] \ . \] 
Again we devide it into three parts 
\[  \mathbb{E}\big[\chi(Y^n_{1/2}); Z_n>0\big] = \sum_{k=0}^{m-1} \ldots + \sum_{k=m}^{n-m} \ldots + \sum_{j=0}^{m-1} \ldots \]
with $m < \lfloor n/2 \rfloor$ fixed. 

For the terms in the first sum $\sum_{k=0}^{m-1} \ldots$ we use the formula
\begin{align*} \mathbb{E}\big[\chi(Y^n_{1/2}); Z_n>0 , \ &\tau_k=k, \min_{k<l\le n} S_l \ge S_k\big] \\ &=  \mathbb{E}\big[\xi_{n,k}(Z_k,S_k);\tau_k=k\big] \mathbb{P}\{L_{n-k}\geq 0] 
\end{align*}
with 
\[
\xi_{n,k}(z,r) = \mathbb{E}\big[ \chi\big(e^{-r} e^{-S_{m_{n-k}}} Z_{m_{n-k}}\big) ;Z_{n-k}>0\ \big |\ L_{n-k}\geq 0, Z_0=z\big] \ . \] 
and $m_n = \lfloor (n+k)/2\rfloor -k$. From Lemma \ref{lem1} we see that this expression is convergent for $n \to \infty$, thus $\sum_{k=0}^{m-1} \ldots $ can be treated just as in the proof of Lemma \ref{lem3}.

For the terms in the second sum $\sum_{k=m}^{n-m} \ldots$ we use the estimate
\begin{align*}\big|\mathbb{E}\big[\chi(Y^n_{1/2}); Z_n>0 , \ &\tau_k=k, \min_{k<l\le n} S_l \ge S_k\big] \big|\\ &\le \sup |\chi| \mathbb{P} \big\{Z_n>0 , \tau_k=k, \min_{k<l\le n} S_l \ge S_k\big\} \ .
\end{align*}
As in Lemma \ref{lem3} we may conclude that the second sum becomes negligible by choosing $m$ sufficiently large.

For the terms of the third sum $\sum_{j=0}^{m-1} \ldots$ we use the formula
\begin{align*} \mathbb{E}\big[\chi(Y^n_{1/2}); Z_n>0 ,\tau_{n-j}=n-j&, \min_{n-j< l \le n} S_l \ge S_{n-j}\big] \\ &= \mathbb{E}\big[\chi(Y^n_{1/2}) \xi_{j}(Z_{n-j});\tau_{n-j}=n-j\big]
\end{align*}
with
\[
\xi_j(z) = \mathbb{P}\{Z_j>0, L_j\ge 0\ |\ Z_0=z\}  \ .
\]
Now we may apply Lemma \ref{lem2}. Altogether \eqref{0505} is proven. Also the statement $W>0$ a.s. follows from Lemma \ref{lem1} and \ref{lem2}.

It remains to show that for $\varepsilon > 0$
\begin{align}
\lim_{n\rightarrow\infty} \mathbb{P}\big\{ \sup_{0 \le t \le 1} |Y_t^n-Y^n_0|>\varepsilon\ \big|\ Z_n>0\big\}= 0 \ . \label{1505}
\end{align}
First we consider a fixed environment. Then $(e^{-S_ i} Z_i)$ is a martingale, thus applying the Doob inequality  we obtain for $k < r_n$
\begin{align*}
\mathbb{P}\big\{ \sup_{0 \le t \le 1} |Y_t^n-Y^n_0|>&\ \varepsilon\ \big|\ Z_k=z,\Pi\big\}  \\ &\leq \ \varepsilon^{-2} \mathbb{E}\big[ (e^{-S_{n-r_n}} Z_{n-r_n} -e^{-S_{r_n}} Z_{r_n})^2\ |\ Z_k=z,\Pi\big] \ .
\end{align*}
Also a straightforward calculation gives for $i \ge k$
\begin{align*} \mathbb{E}\big[ \big(e^{-S_{i+1}} Z_{i+1}-&\ e^{-S_i} Z_i\big)^2\ \big|\ Z_k=z, \Pi\big]  \\ &=   ze^{-S_k}(\eta_{i+1} e^{-S_i} + e^{-S_{i+1}} -e^{-S_{i}}) \ .
\end{align*}
Given $\Pi$ and $Z_k=z$ the process $(e^{-S_i} Z_i)$ is therefore a $L_2$-martingale. Consequently
\begin{align*}
\mathbb{E}\big[ (e^{-S_{n-r_n}} Z_{n-r_n} &-e^{-S_{r_n}} Z_{r_n})^2\ \big|\ Z_k=z, \Pi\big] \\ &= \sum_{i=r_n}^{n-r_n-1} \mathbb{E}\big[ \big(e^{-S_{i+1}} Z_{i+1}- e^{-S_i} Z_i\big)^2\ \big|\ Z_k=z,\Pi\big] \\ &\le ze^{-S_k}\Big(\sum_{i=r_n}^{n-r_n-1} \eta_{i+1} e^{-S_i} + e^{-S_{n-r_n}} \Big) \ .
\end{align*}
Letting
\begin{align*}
&U_n := \sum_{i=r_n}^{\lfloor n/2 \rfloor} \eta_{i+1} e^{-S_i} \ , \ V_n := \sum_{i=r_n+1}^{\lceil n/2 \rceil} \eta_{i} e^{S_i}+ e^{S_{r_n}}\\
&\varphi(u,v,x) :=  \varepsilon^{-2}(u+ v e^{-x})^+ \ ,
\end{align*}
we obtain altogether (recall the definition of $\tilde V_n$ from the proof of Lemma \ref{lem1})
\[ \mathbb{P}\big\{ \sup_{0 \le t \le 1} |Y_t^n-Y^n_0|> \varepsilon \ \big|\ Z_k=z,\Pi\big\}  \le \big(ze^{-S_k}\varphi(U_n, \tilde V_n,S_n)\big)\wedge 1 \ . \]

Once again we proceed in the by now established manner:
\begin{align*}
\mathbb{P}\big\{ &\sup_{0 \le t \le 1} |Y_t^n-Y^n_0|> \varepsilon ,\ Z_n >0\}\\ & = \sum_{k=0}^n \mathbb{P}\big\{ \sup_{0 \le t \le 1} |Y_t^n-Y^n_0|> \varepsilon ,\ Z_n>0 , \tau_k=k, \min_{k<l\le n} S_l \ge S_k\big\}\\ &\le \sum_{k=0}^{m-1} \mathbb P\big\{ \sup_{0 \le t \le 1} |Y_t^n-Y^n_0|> \varepsilon  ,\ \tau_k=k, \min_{k<l\le n} S_l \ge S_k\big\}
\\ & \quad \mbox{} + \sum_{k=m}^{n-m}\mathbb{P}\big\{Z_n>0 , \tau_k=k, \min_{k<l\le n} S_l \ge S_k\big\}
\\ & \quad \mbox{} + \sum_{j=0}^{m-1} \mathbb{P}\big\{ \sup_{0 \le t \le 1} |Y_t^n-Y^n_0|> \varepsilon ,\ Z_{n-j}>0 , \tau_{n-j}=n-j \big\} \ .
\end{align*}
As to the sum $\sum_{k=0}^{m-1} \ldots$ 
\begin{align*}
\mathbb P\big\{  \sup_{0 \le t \le 1}& |Y_t^n-Y^n_0|> \ \varepsilon ,\tau_k=k, \min_{k<l\le n} S_l \ge S_k\big\} \\ &\le \mathbb E\big[\mathbb{P}\big\{ \sup_{0 \le t \le 1} |Y_t^n-Y^n_0|> \varepsilon \ \big|\ Z_k,\Pi\big\}; \min_{k<l\le n} S_l \ge S_k \big]\\ &\le \mathbb E\big[ \big(Z_ke^{-S_k}\varphi(U_n, \tilde V_n,S_n)\big)\wedge 1 ;\min_{k<l\le n} S_l \ge S_k \big]  \ . 
\end{align*}
Next $\mathbb E[(Z_ke^{-S_k}\varphi(U_n, \tilde V_n,S_n))\wedge 1\, |\, \Pi]\le \varphi(U_n, \tilde V_n,S_n)\wedge 1$ because of Jensen's inequality applied to the concave function $x \mapsto x\wedge 1$. Thus we obtain
\begin{align*}
\mathbb P\big\{  \sup_{0 \le t \le 1} |Y_t^n&-Y^n_0|>  \varepsilon ,\tau_k=k, \min_{k<l\le n} S_l \ge S_k\big\} \\ 
&\le \mathbb E[ \varphi(U_{n}, \tilde V_{n},S_{n})\wedge 1; \min_{k<l\le n} S_l \ge S_k ]\\
&= \mathbb E\big[ \big(e^{-S'} \varphi(U_{n-k}, \tilde V_{n-k},S_{n-k})\big)\wedge 1; L_{n-k} \ge 0 \big] \ ,
\end{align*} 
where in the last expectation $S'$ is distributed as $S_k$ before and independent from the other terms and also $r_{n-k}$ is replaced by $r_n-k$. Now $S_n \to -\infty$ $\mathbf P^-$-a.s. (compare Lemma 2.6 in \cite{agkv}). This together with Lemma \ref{le2} gives that $U_n \to 0 $ $\mathbf P^+$-a.s. and $V_n \to 0$ $\mathbf P_x^-$-a.s. From Theorem \ref{pr4} it follows that
\[\mathbb P\big\{  \sup_{0 \le t \le 1} |Y_t^n-Y^n_0|> \ \varepsilon ,\tau_k=k, \min_{k<l\le n} S_l \ge S_k\big\}  = o\big( \mathbb P \{ L_n \ge 0 \} \big)\ . \]
Thus the sum $\sum_{k=0}^{m-1} \ldots$ is negligible. As we already know, the same is true for $\sum_{k=m}^{n-m} \ldots$ by choosing $m$ large. As to $\sum_{j=0}^{m-1} \ldots $ we obtain by means of H\"older's inequality for conjugate numbers $p,q > 1$
\begin{align*}
\mathbb{P}\big\{ \sup_{0 \le t \le 1} &|Y_t^n-Y^n_0|>\varepsilon , \  Z_n>0 \ \big| \  \Pi \big\} \\ &\le \mathbb{E}\big[Z_n^{1/p};\sup_{0 \le t \le 1} |Y_t^n-Y^n_0|>\varepsilon \ \big| \  \Pi \big]\\
& \le  
\mathbb E[Z_n | \Pi]^{1/p}  \mathbb{P}\big\{ \sup_{0 \le t \le 1} |Y_t^n-Y^n_0|>\varepsilon\  \ \big| \  \Pi \big\}^{1/q}  \\
& \le \exp\big( S_n/p) (\varphi(U_n, \tilde V_n,S_n)\wedge 1)^{1/q}  
\end{align*}
and consequently
\begin{align*}
\mathbb{P}\big\{  \sup_{0 \le t \le 1} &|Y_t^n-Y^n_0|> \varepsilon ,\ Z_{n-j}>0 , \tau_{n-j}=n-j \big\} \\ &\le \mathbb E\big[ \exp\big( S_{n-j}/p) (\varphi(U_{n-j}, \tilde V_{n-j},S_{n-j})\wedge 1)^{1/q} ; \tau_{n-j}=n-j \big] \ ,
\end{align*}
where $r_{n-j}$ is again replaced by $r_n-j$. 
Now we choose $p$ such that $1/p > \beta$. Then, switching once more to the tilted measure $\mathbf E$, we may apply Theorem \ref{pr5} again to obtain
\[ \mathbb{P}\big\{  \sup_{0 \le t \le 1} |Y_t^n-Y^n_0|> \varepsilon ,\ Z_{n-j}>0 , \tau_{n-j}=n-j \big\} = o\big( \mathbb P \{ L_n \ge 0 \} \big)\ . \]
Altogether
\[  \mathbb{P}\big\{ \sup_{0 \le t \le 1} |Y_t^n-Y^n_0|>\varepsilon\ \big|\ Z_n>0\big\} = o\big( \mathbb P \{ L_n \ge 0 \} \big) = o(1) \]
and this gives \eqref{1505}.
\qed

\bigskip      
Department of       
Discrete Mathematics \\ Steklov Institute \\ 8 Gubkin Street,      
117\,966 Moscow, GSP-1, Russia \\ viafan@mail.ru, vatutin@mi.ras.ru    \\

\noindent
Fachbereich Mathematik\\ Universit\"at Frankfurt\\ Fach   
    187, D-60054 Frankfurt am Main, Germany \\         
 kersting@math.uni-frankfurt.de , boeinghoff@math.uni-frankfurt.de

  \end{document}